\documentclass[reqno,english]{amsart}
\usepackage{amsfonts,amsmath,latexsym,verbatim,amscd,mathrsfs,color,array}
\usepackage[colorlinks=true]{hyperref}
\usepackage{amsmath,amssymb,amsthm,amsfonts,graphicx,color}
\usepackage{amssymb}
\usepackage{pdfsync}
\usepackage{epstopdf}
\usepackage{cite}
\usepackage{graphicx}
\usepackage{datetime}

\usepackage{setspace}

\usepackage{caption2}
\usepackage{subfigure}
\usepackage{float}

\newcommand{\bt}{\begin{theorem}}
\newcommand{\et}{\end{theorem}}
\newcommand{\bl}{\begin{lemma}}
\newcommand{\el}{\end{lemma}}
\newcommand{\bd}{\begin{definition}}
\newcommand{\ed}{\end{definition}}
\newcommand{\bc}{\begin{corollary}}
\newcommand{\ec}{\end{corollary}}
\newcommand{\bp}{\begin{proof}}
\newcommand{\ep}{\end{proof}}
\newcommand{\bx}{\begin{example}}
\newcommand{\ex}{\end{example}}
\newcommand{\bi}{\begin{exercise}}
\newcommand{\ei}{\end{exercise}}
\newcommand{\bo}{\begin{prop}}
\newcommand{\eo}{\end{prop}}
\newcommand{\br}{\begin{remark}}
\newcommand{\er}{\end{remark}}
\newcommand{\be}{\begin{equation}}
\newcommand{\ee}{\end{equation}}
\newcommand{\ba}{\begin{align}}
\newcommand{\ea}{\end{align}}
\newcommand{\bn}{\begin{enumerate}}
\newcommand{\en}{\end{enumerate}}
\newcommand{\bg}{\begin{align*}}
\newcommand{\bcs}{\begin{cases}}
\newcommand{\ecs}{\end{cases}}

\newcommand{\bean}{\begin{eqnarray*}}
\newcommand{\eean}{\end{eqnarray*}}

\newtheorem{example}{Example}[section]
\newtheorem{definition}{Definition}[section]
\newtheorem{theorem}{Theorem}[section]
\newtheorem{lemma}{Lemma}[section]

\newtheorem{prop}{Proposition}[section]
\newtheorem{remark}{Remark}[section]

\numberwithin{equation}{section}

\begin{document}
\title[Symplectic algorithms]{Symplectic algorithms for stable manifolds in control theory}
\author[G. Chen]{Guoyuan Chen}
\address{\noindent
School of Data Sciences, Zhejiang University of Finance \& Economics, Hangzhou 310018, Zhejiang, P. R. China}
\email{gychen@zufe.edu.cn}

\author[G. Zhu]{Gaosheng Zhu}
\address{\noindent
School of Mathematics, Tianjin University, Tiajin 300072, P. R. China}
\email{gaozsc@163.com}

\thanks{The authors are grateful to the anonymous referees for useful comments and suggestions. G. Chen is supported by NSFC (No.11771386) and First Class Discipline of Zhejiang - A (Zhejiang University of Finance and Economics- Statistics). G. Zhu is supported by NSFC (No.11871356).}
\begin{abstract}
In this note,  we propose a symplectic algorithm for the stable manifolds of the Hamilton-Jacobi equations combined with an iterative procedure in [Sakamoto-van~der Schaft,
IEEE Transactions on Automatic Control, 2008]. Our algorithm includes two key aspects. The first one is to prove a precise estimate for radius of convergence and the errors of local approximate stable manifolds. The second one is to extend the local approximate stable manifolds to larger ones by symplectic algorithms which have better long-time behaviors than general-purpose schemes. Our approach avoids the case of divergence of the iterative sequence of approximate stable manifolds, and reduces the computation cost. We illustrate the effectiveness of the algorithm by an optimal control problem with exponential nonlinearity.
\end{abstract}
\maketitle


\section{Introduction}

It is well known that an optimal feedback control can be given by solving an associated Hamilton-Jacobi (HJ) equation (see e.g. \cite{lee1967foundations}) and $H^\infty$ feedback control can be obtained from solutions of one or two Hamilton-Jacobi equations (see e.g. \cite{van1991state, isidori1992disturbance, van19922, ball1993h}).
Unfortunately, the Hamilton-Jacobi equation in general can not be solved analytically. Hence numerical method becomes important.
Seeking approximate solutions of Hamilton-Jacobi equations from control theory has been studied extensively. There are several approaches: Taylor series method, Galerkin method, state-dependent Riccati equation method, algebraic method, etc. See e.g. \cite{al1961optimal,lukes1969optimal, kreisselmeier1994numerical,beard1997galerkin,bea1998successive, mracek1998control, markman2000iterative, beeler2000feedback,aliyu2003approach, aliyu2003transformation, mceneaney2007curse, navasca2007patchy, hunt2010improved, ohtsuka2010solutions, aguilar2014numerical} and the references therein. These methods may have good performance for concrete control systems. However, in general, they may have various disadvantages such as heavy computation cost for higher dimensional state spaces, restriction on simple nonlinearity of the systems, etc.

For the stationary Hamilton-Jacobi equations which are related to infinite horizon optimal control and $H^\infty$ control problems, \cite{sakamoto2008} developed an iterative procedure to construct an approximate sequence that converges to the exact solution of the associated Hamiltonian system on the stable manifold. It is based on the fact that the stabilizing solutions of stationary Hamilton-Jacobi equations correspond to the generating functions of the stable manifolds (Lagrangian) of the associated Hamiltonian systems at certain equilibriums (cf. e.g. \cite{markman2000iterative,sakamoto2002analysis, sakamoto2008}). This approach has better performances for various nonlinear feedback control systems, especially for the ones with more complicated nonlinearities, see e.g. \cite{sakamoto2013case,horibe2017optimal, horibe2018nonlinear}.

We should note that the computation approach in \cite{sakamoto2008} (as well as \cite{sakamoto2013case, horibe2017optimal, horibe2018nonlinear}) depends essentially on the radius of convergence of the iterative procedure which is not estimated analytically. Moreover, since the errors of less iterative steps are tremendous when we enlarge the local approximate stable manifolds in unstable direction, to obtain a stable manifold with proper size for applications, the number of iterative steps should be large. This may make the computation time-consuming.

In this note, inspired by \cite{sakamoto2008}, we  construct a sequence of local approximate stable manifolds of stationary Hamilton-Jacobi equations by iterative procedure with a precise estimate of radius of convergence, and enlarge the local manifolds to larger ones by symplectic algorithms. To be more precise, as in \cite{sakamoto2008} we firstly generate a sequence of local approximate stable manifolds near the equilibrium by an iterative procedure to solve the associated Hamiltonian system of the Hamilton-Jacobi equation. Then we extend the local approximate stable manifolds to large ones by symplectic algorithms for solving initial value problem of the associated Hamiltonian system (Section \ref{s:symplectic} below). We emphasize that in our approach, the radius of convergence of the iterative sequences are estimated precisely and the error of the local approximate stable manifold can be controlled as small as possible (Theorem \ref{t:iteration} below).
Therefore the significant point is how to extend the local stable manifolds.
There are various numerical methods for general ODEs, e.g. Runge-Kutta methods of various orders. For our applications in Hamiltonian systems, we will use symplectic algorithms.
The symplectic structure plays a significant role in design of numerical methods for Hamiltonian systems. Generally speaking, the symplectic algorithms are designed to preserve the symplectic structure of the Hamiltonian systems. Hence, for Hamiltonian systems, symplectic algorithm improves qualitative behaviours, and gives a more accurate long-time integration comparing with general-purpose methods such as Runge-Kutta schemes. Various kinds of symplectic algorithms, e.g., symplectic Euler, St\"ormer-Verlet, symplectic Runge-Kutta, were constructed since 1950s. A detailed history of symplectic methods and related topics can be found in \cite{gauckler2017dynamics}.
We refer the readers to the books \cite{hairer2006geometric, feng2010symplectic, blanes2016concise, sanz2018numerical} for a complete presentation of various symplectic algorithms for Hamiltonian systems.

The note is organized as follows. In Section \ref{s:Ham-Joc-control}, basic notations for the stable manifolds of Hamilton-Jacobi equations in control theory are given. Section \ref{s:iteration} is devoted to constructing the iterative procedure, and proves the precise estimate of radius of convergence as well as the error of the approximate solutions. The symplectic scheme which extends the local approximate stable manifold is described in Section \ref{s:symplectic}. In Section \ref{s:example}, the effectiveness of our algorithm is illustrated by application to an optimal control problem. The Appendix gives the details of the proof of Theorem \ref{t:iteration}.

\section{The Hamilton-Jacobi equation and the stable manifolds}\label{s:Ham-Joc-control}
In this section, for the convenience of the reader we recall the relevant material without proofs, thus making our exposition self-contained. For more details, see \cite{sakamoto2008} and \cite{sakamoto2013case}.

Let $\Omega\subset\mathbb R^d$ be a domain containing $0$. We consider the following Hamilton-Jacobi equation
\begin{eqnarray}\label{e:HJ}
H(x,p)=p^Tf(x)-\frac{1}{2}p^TR(x)p+q(x)=0,
\end{eqnarray}
where $p := \nabla V$ for some unknown function $V$, $f:\Omega\to \mathbb R^d$, $R:\Omega\to \mathbb R^{d\times d}$, $q:\Omega\to \mathbb R$ are $C^\infty$ and $R(x)$ is symmetric matrix for all $x\in\Omega$. Furthermore, we assume that $f(0)=0$ and $q(0)=0$, $\frac{\partial q}{\partial x}(0)=0$. Hence for $x$ near $0$,
$
f(x)=Ax+O(|x|^2),\,q(x)=\frac{1}{2}x^TQx+O(|x|^3),
$
where $A=\frac{\partial f}{\partial x}(0)$, $Q=\frac{\partial^2q}{\partial x^2}(0)$ is the Hessian of $q$ at $0$.
We say that a solution $V$ of \eqref{e:HJ} is \emph{stabilizing} if $p(0)=0$ and $0$ is an asymptotically stable equilibrium of the vector field $f(x)-R(x)p(x)$ where $p(x)=\nabla V(x)$.

This kind of problem arises from infinite horizon optimal control. For example, consider the optimal control problem
\begin{eqnarray}
\dot{x}=f(x)+g(x)u,\quad \mbox{in }\Omega,
\end{eqnarray}
where $g$ is a smooth $d\times m$ matrix-valued function, $u$ is an $m$-dimensional feedback control function of column form.
Assume the instantaneous cost function
$
L(x,u)=\frac{1}{2}(q(x)+ru^Tu),
$
where $r>0$, $q(x)\ge 0$ is a smooth function. Define the cost functional by
$
J(x,u)= \int_0^{+\infty}L(x(t),u(t))dt.
$
Then the corresponding Hamilton-Jacobi equation is of form \eqref{e:HJ}
with $R(x):=r^{-1}g(x)g(x)^T$.
The optimal controller is
\begin{eqnarray}\label{e:u*12}
u(x)=-r^{-1}g(x)^T p(x),
\end{eqnarray}
where $p(x)=\nabla V(x)$.

From the symplectic geometry point of view, a stabilizing solution $V$ of \eqref{e:HJ} corresponds to a stable (Lagrangian) submanifold. That is,
$
\Lambda_V:=\{(x,p)\,|\,p=\nabla V\}
$
is a stable (Lagrangian) manifold which is invariant under the flow of the associated Hamiltonian system of \eqref{e:HJ}:
\begin{eqnarray}\label{e:ham-flow}
\left\{\begin{array}{l}
         \dot x=H_p(x,p)=f(x)-R(x)p, \\
         \dot p=-H_x(x,p)=-(\frac{\partial f}{\partial x})^Tp + \frac{1}{2}\frac{\partial(p^TR(x)p)}{\partial x}-(\frac{\partial q}{\partial x})^T.
       \end{array}
\right.
\end{eqnarray}
Conversely, if a $d$-dimensional manifold $\Lambda$ in $(x,p)$-space is invariant with respect to the flow \eqref{e:ham-flow}, and at some point $(x_0,p_0)$, the projection of $\Lambda$ to the $x$-space is surjective, then $\Lambda$ is a Lagrangian submanifold in a neighborhood of $(x_0,p_0)$ and there is a solution $V$ of \eqref{e:HJ} in a neighborhood of $x_0$ such that $\Lambda_V=\Lambda$. Therefore, if we get the stable manifold, then the optimal controller can be obtained from \eqref{e:u*12}. See e.g. \cite{sakamoto2008}.

Denote $z=(x,p)$. Let $J$ be the standard symplectic matrix $\left[\begin{array}{cc}
                                                                                                   0 & I_d \\
                                                                                                   -I_d & 0
                                                                                                 \end{array}
\right]$, where $I_d$ denotes the identity matrix of $d$-dimensional. Then the vector field on left side of \eqref{e:ham-flow} is $J^{-1}\nabla H(z)$. Let $X_H(z)=J^{-1}\nabla H(z)$ be the Hamiltonian vector field of $H$. Note that $z_0=(0,0)$ is an equilibrium of $X_H$. Then the derivative of the Hamiltonian vector field at $z_0$ is a Hamiltonian matrix, i.e., $(JDX_H(z_0))^T=JDX_H(z_0)$.
We say that the equilibrium $z_0$ is \emph{hyperbolic} if $DX_H(z_0)$ has no imaginary eigenvalues.
It is well known that if $DX_H(z_0)$ is hyperbolic, then its eigenvalues are symmetric with respect to the imaginary axis. From the Stable Manifold Theorem, there exists a global stable manifold $\mathcal S_{z_0}$ of $z_0$. Moreover, $\mathcal S_{z_0}$ is a Lagrangian submanifold of $(\mathbb R^{2d},\omega)$ where $\omega$ is the standard symplectic structure (see e.g. \cite{AbrMar1978}).
Near $z_0$, the Hamiltonian system \eqref{e:ham-flow} can be rewritten as
\begin{eqnarray}\label{e:local}
\dot z = DX_H(z_0)z+N(z),
\end{eqnarray}
where $N(z)$ is the nonlinear term.

A sufficient condition for the existence local stabilizing solution for \eqref{e:HJ} is obtained by van der Schaft \cite{van1991state} based on an observation on the Riccati equation. Assume $(x_0,p_0)=(0,0)$ without loss of generality. Let $R:=R(0)$. The Riccati equation
\begin{eqnarray}\label{e:Riccati}
PA+A^TP-PRP+Q=0
\end{eqnarray}
is the linearization of \eqref{e:HJ} at the origin. A symmetric matrix $P$ is said to be the \emph{stabilizing} solution of \eqref{e:Riccati} if it is a solution of  \eqref{e:Riccati} and $A-RP$ is stable.
The Riccati equation \eqref{e:Riccati} has a stabilizing solution if and only if the following two conditions hold:
(1) $DX_H(z_0)$ is hyperbolic; (2) the generalized eigenspace $E_-$ for $d$ stable eigenvalues satisfies
  $
  E_-\oplus {\rm Im}(0,I_d)^T=\mathbb R^{2d}.
  $
See e.g. \cite{lancaster1995algebraic}. If \eqref{e:Riccati} has a stabilizing solution $P$, then there exists a local stabilizing solution $V$ of \eqref{e:HJ} around the origin such that $\frac{\partial^2V}{\partial x^2}(0)=P$. This yields a local solution of Hamilton-Jacobi solution (\cite{sakamoto2008}).
Consider the Lyapunov equation
\begin{eqnarray}\label{e:Lyapunov1}
(A-RP)S+S(A-RP)^T=R.
\end{eqnarray}
Some direct computations yield the following result (\cite{sakamoto2008}).
\begin{lemma}\label{l:diag}
Assume that $S$ is a solution of \eqref{e:Lyapunov1}. Then
\begin{eqnarray}
T^{-1}DX_H(z_0)T=\left[
                   \begin{array}{cc}
                     B & 0 \\
                     0 & -B^T \\
                   \end{array}
                 \right],
\end{eqnarray}
where $B=A-RP$ and
$
T=\left[
    \begin{array}{cc}
      I_d & S \\
      P & PS+I_d \\
    \end{array}
  \right].
$
\end{lemma}

From Lemma \ref{l:diag}, by coordinates transformation
$\left[\begin{array}{c}
             \bar x \\
            \bar p \\
           \end{array}
         \right]=T^{-1}\left[
           \begin{array}{c}
             x \\
             p \\
           \end{array}
         \right]$,
system \eqref{e:local} becomes a separated form
\begin{eqnarray}\label{e:Ham-transform}
\left\{\begin{array}{l}
         \dot {\bar x} = B\bar x + N_s(\bar x,\bar p) \\
         \dot {\bar p} = -B^T\bar p+N_u(\bar x,\bar p)
       \end{array}\right.,
\end{eqnarray}
where $N_s(\bar x,\bar p)$ and $N_u(\bar x,\bar p)$ are nonlinear terms corresponding to $N(z)$ after the coordinates transformation.

\section{The local approximate stable manifolds: iteration}\label{s:iteration}

In this section, we shall give an iterative procedure to construct a sequence of local approximate stable manifolds near equilibrium for the Hamiltonian system in the form
\begin{eqnarray}\label{e:after-ct1}
\left\{\begin{array}{l}
         \dot x = Bx + n_s(x,y),\\
         \dot y = -B^Ty+n_u(x,y).
       \end{array}\right.
\end{eqnarray}

{\bf Assumption 1:} $B$ has eigenvalues with negative real parts. It follows that there exist positive constants $a, b$ such that $\|e^{Bt}\|\le a e^{-bt}$ for $t\ge 0$.

{\bf Assumption 2:} $n_s,n_u:\mathbb R^d\times\mathbb R^d\to \mathbb R^d$ are continuous and satisfy the following conditions:
For all $L>0$, $0<l\le L$, $|x|+|y|\le l$ and $|x'|+|y'|\le l$,
  \begin{eqnarray}\label{e:n-assumption}
  &&|n_s(x,y)-n_s(x',y')|\le M(L) l(|x-x'|+|y-y'|),\notag\\
  &&|n_u(x,y)-n_u(x',y')|\le M(L) l(|x-x'|+|y-y'|),\notag
  \end{eqnarray}
where $M(L)$ is increasing with respect to $L$.

As in \cite{sakamoto2008}, to solve equation \eqref{e:after-ct1},
we define the following iterative sequence for $k=1,2,\cdots$,
\begin{eqnarray}\label{e:iterative-int}
\left\{\begin{array}{l}
         x_{k+1}=e^{Bt}\xi+\int_0^t e^{B(t-s)}n_s(x_k(s),y_k(s))ds\\
         y_{k+1}=-\int_t^{\infty}e^{B^T(t-s)} n_u(x_k(s),y_k(s))ds
       \end{array}
\right.
\end{eqnarray}
and $x_0 = e^{Bt}\xi$, $y_0=0$ for an arbitrary $\xi\in \mathbb R^d$. Equivalent to \eqref{e:iterative-int}, we consider the following ODE:
\begin{eqnarray}\label{e:iterative}
\left\{\begin{array}{l}
         \dot x_{k+1}=B x_{k+1}+n_s(x_k(t),y_k(t))\\
         \dot y_{k+1}=-B^Ty_{k+1}+ n_u(x_k(t),y_k(t))
       \end{array}
\right.
\end{eqnarray}
with boundary conditions $x_{k+1}(0)=\xi,\,y_{k+1}(+\infty)=0$ and $x_0=e^{Bt}\xi$, $y_0=0$, $t\ge 0$. This form is more convenient to apply numerical methods to ODEs.

Inspired by \cite[Theorem 5]{sakamoto2008}, we can prove the following convergence result whose proof is included in Appendix.
\begin{theorem}\label{t:iteration}
Assume that system \eqref{e:after-ct1} satisfies Assumption 1-2, $M(L)L\ge\frac{3b}{8a}$ and $|\xi|\le \frac{3b}{16a^2M}$, where $M=M(L)$ is the constant depending on $L$ given by Assumption 2. Let $\{x_k\}$ and $\{y_k\}$ be the sequences defined by \eqref{e:iterative-int}. Then $x_k(t)\to 0$, $y_k(t)\to 0$ as $t\to +\infty$, and there exist functions $x$ and $y$ such that $\{x_k\}$ and $\{y_k\}$ uniformly converge to $x$ and $y$ respectively, $(x,y)$ is solution of \eqref{e:after-ct1},  and for all $t\ge 0$,
\begin{eqnarray}\label{e:x-y}
&&|x_{k}(t)-x(t)|
\le C(a,b,M)\frac{|\xi|^2}{2^{k-1}}e^{-bt},\notag\\
&&|y_{k}(t)-y(t)|\le C(a,b,M)\frac{|\xi|^2}{2^{k-1}}e^{-2bt},
\end{eqnarray}
where
$
C(a,b,M)>0
$ is a constant depending only on $a,b,M$.
\end{theorem}
\begin{remark}\label{r:r-error}
Compared to \cite[Theorem 5]{sakamoto2008}, we improve the result in two aspects: (1) a sufficient estimate of $|\xi|$ is given; (2) the error of iteration is estimated precisely. The constant $C(a,b, M)$ also can be calculated explicitly (see the proof of Theorem \ref{t:iteration} in the Appendix below).
\end{remark}

\section{Extension of the local stable manifolds by symplectic algorithms}\label{s:symplectic}

In this section, the local stable manifolds will be enlarged by symplectic algorithms.

\subsection{Structure of the approximate stable manifold}
By Lemma \ref{l:diag} and Theorem \ref{t:iteration}, we obtain a sequence of local approximate stable manifolds of \eqref{e:ham-flow} near equilibrium $(0,0)$. Let
$
\mathbb S_{\rho}=\{\xi\in \mathbb R^d\,|\,|\xi|=\rho\},
$
where $\rho$ is the radius of convergence chosen by Theorem \ref{t:iteration}.
Denote the local approximate stable manifold by
$$
\Lambda_k=\{(x_k(t,\xi),p_k(t,\xi))\,|\,t\ge 0,\,\xi\in \mathbb S_\rho\},
$$
and denote the boundary of $\Lambda_k$ by $\partial \Lambda_k$. Then $\partial \Lambda_k=\{(\xi,p_k(0,\xi))\,|\,\xi\in \mathbb S_\rho\}$.
As $k\to \infty$, $\Lambda_k$ tends to the exact stable manifold near equilibrium $(0,0)$ parameterized by $(t,\xi)$
\begin{eqnarray}
\Lambda:=\{(x(t,\xi),p(t,\xi))\,|\,t\ge 0,\,\xi\in \mathbb S_\rho\}.
\end{eqnarray}
Its boundary $\partial \Lambda=\{(\xi,p(0,\xi))\,|\,\xi\in \mathbb S_\rho\}$.
Consider the initial value problem of Hamiltonian system
\begin{eqnarray}\label{e:after-ct-extension}
\left\{\begin{array}{l}
         \dot x = H_p(x,p) \\
         \dot p = -H_x(x,p)
       \end{array}\right.,
       \mbox{ for } t\le 0\mbox{ with } (x(0),p(0))\in \partial\Lambda.
\end{eqnarray}
Then by the invariance of the stable manifold,
\begin{eqnarray}
\Lambda_{\rm g}&:=&\{(x(t),p(t))\,|\,t\le 0,\,(x(0),p(0))\in \partial \Lambda\}\cup \Lambda \notag\\
&=&\{(x(t,\xi),p(t,\xi))\,|\,t\in \mathbb R,\,\xi\in \mathbb S_\rho\},\notag
\end{eqnarray}
is the global stable manifold of $(0,0)$. Hence we extend local stable manifold $\Lambda$ to the global one. Moreover, we should emphasize that each solution curve of \eqref{e:after-ct-extension}, $(x(t),p(t))$, with $(x(0),p(0))\in \partial \Lambda$ lies in $\Lambda_{\rm g}$.
Numerically, we compute the approximations by the following initial value problem
\begin{eqnarray}\label{e:after-ct-extension-app-H}
\left\{\begin{array}{l}
         \dot x = H_p(x,p) \\
         \dot p = -H_x(x,p)
       \end{array}\right.
       \mbox{for } t\le 0\mbox{ with } (x_k(0),p_k(0))\in \partial \Lambda_k.
\end{eqnarray}
Letting $(x_k(t),p_k(t))$ be numerical solutions of \eqref{e:after-ct-extension-app-H}, we obtain an approximate stable manifold
\begin{eqnarray}\label{e:app-stable}
\Lambda_{k,\rm g}:=\{(x_k(t),p_k(t))\,|\,t\le 0,\,(x_k(0),p_k(0))\in \partial \Lambda_k\}\cup \Lambda_k\notag\\
=\{(x_k(t,\xi),p_k(t,\xi))\,|\,t\in \mathbb R,\,\xi\in \mathbb S_\rho\}.
\end{eqnarray}
Therefore, the key point is to numerically solve the problem \eqref{e:after-ct-extension-app-H}. For general ODEs, there are many kinds of numerical algorithms for the initial value problems. For example, Runge-Kutta of various orders. However, we should point out that \eqref{e:after-ct-extension-app-H} is a Hamiltonian system. Symplectic algorithms have better performance for such kind of systems.

\subsection{Symplectic algorithms of Hamiltonian systems}
For Hamiltonian systems, it is natural to use symplectic algorithms. A numerical one-step method $y_{n+1}=\Phi_h(y_n)$ (with step size $h$) is called symplectic if $\Phi_h$ is a symplectic map, that is,
$
D\Phi_h(y)^TJD\Phi_h(y)=J,
$
where $D\Phi_h(y)$ is the tangent map of $\Phi_h$ at $y$. Symplectic structure is an essential property of Hamiltonian system (cf., e.g., \cite[Chapter VI.2]{hairer2006geometric}). Symplectic algorithms preserve this geometric structure for each step. Hence compared to general purpose numerical algorithms (e.g. Runge-Kutta methods), symplectic algorithms have much better long-time qualitative behaviours.
There are many types of symplectic algorithms, e.g., symplectic Runge-Kutta of various orders, St\"ormer-Verlet methods, Nystr\"om method, etc. We refer the readers to \cite{hairer2006geometric,feng2010symplectic} for more details of symplectic algorithms.

In what follows, we illustrate our procedure of extension of the local stable manifold by the St\"ormer-Verlet method which is a simple symplectic algorithm of 2-order. Other symplectic algorithms of higher orders may have better numerical results.

Let $h$ be a step size, and $t_0=0$ be the initial time, and let $t_n=nh$, $t_{n+1/2}=(n+1/2)h$. Hence for problem \eqref{e:after-ct-extension-app-H}, we should choose $h<0$. Let $p_{n}=p(t_n)$, $p_{n+1/2}=p(t_{n+1/2})$,
$x_{n}=x(t_n)$, $x_{n+1/2}=x(t_{n+1/2})$. Denote $H_p(x,p)=\frac{\partial H}{\partial p}(x,p)$, $H_x(x,p)=\frac{\partial H}{\partial x}(x,p)$. In our case, the Hamiltonian $H(x,p)$ is given by \eqref{e:HJ} where $H_p$ and $H_x$ can be calculated.
Then we have the following theorem.
\begin{theorem}\label{t:SV}
Given $(x_n,p_n)$, the St\"{o}rmer-Verlet schemes
\begin{equation}\label{e:SV1}
\left\{\begin{array}{l}
         p_{n+1/2}=p_n-\frac{h}{2}H_x(x_n,p_{n+1/2}) \\
         x_{n+1} = x_n+\frac{h}{2}\left[H_p(x_n, p_{n+1/2})+H_p(x_{n+1}, p_{n+1/2})\right] \\
         p_{n+1}=p_{n+1/2}-\frac{h}{2}H_x(x_{n+1}, p_{n+1/2}),\quad\mbox{ and }
       \end{array}
\right.
\end{equation}
\begin{equation}\label{e:SV2}
\left\{\begin{array}{l}
         x_{n+1/2}=x_n+\frac{h}{2}H_p(x_{n+1/2}, p_{n}) \\
         p_{n+1} = p_n-\frac{h}{2}\left[H_x(x_{n+1/2}, p_{n})+H_x(x_{n+1/2}, p_{n+1})\right] \\
         x_{n+1}=x_{n+1/2}+\frac{h}{2}H_p(x_{n+1/2}, p_{n+1}),
       \end{array}
\right.
\end{equation}
are symplectic methods of order 2.
\end{theorem}
A complete proof of this Theorem and more details of the St\"ormer-Verlet method can be found in \cite{hairer2003geometric, hairer2006geometric}.
Note that in general \eqref{e:SV1} and \eqref{e:SV2} are implicit equations which can be solved by Newton's iteration method. For example, in \eqref{e:SV1}, the first two equations are implicit. The third equation is explicit if $p_{n+1/2}, x_{n+1}$ were found. Recall that the key point of Newton's iteration method is to give a proper initial guess at the beginning of iteration. For the first two equations of \eqref{e:SV1}, a good initial guess of $(x_{n+1}, p_{n+1/2})$ is $(x_{n}, p_{n})$ since usually the step size $h$ is small and $(x_{n+1}, p_{n+1/2})$ is close to $(x_{n}, p_{n})$. Hence, only a few times of iteration in Newton's method should be applied to arrive at the accuracy needed. Therefore the computation cost is cheap for Newton's iteration at this point.

For our case, from \eqref{e:ham-flow}, the St\"ormer-Verlet scheme \eqref{e:SV1} becomes
\begin{eqnarray}\label{e:sv-s}
\left\{\begin{array}{l}
         p_{n+1/2}=p_n+\frac{h}{2}\left[-\left[\frac{\partial f}{\partial x}(x_n)\right]^Tp_{n+1/2}\right.\\
         \quad\quad\quad\quad+\left.\frac{1}{2}\frac{\partial(p_{n+1/2}^TRp_{n+1/2})}{\partial x}(x_n)-\frac{\partial q}{\partial x}(x_n)\right] \\
         x_{n+1} = x_n+\frac{h}{2}\left[f(x_n)+f(x_{n+1})\right.\\
        \quad\quad\quad\quad \left.-(R(x_n)+R(x_{n+1}))p_{n+1/2}\right] \\
         p_{n+1}=p_{n+1/2}+\frac{h}{2}\left[-\left[\frac{\partial f}{\partial x}(x_{n+1})\right]^Tp_{n+1/2}\right.\\
         \quad\quad\quad\quad+\left.\frac{1}{2}\frac{\partial(p_{n+1/2}^TRp_{n+1/2})}{\partial x}(x_{n+1})-\frac{\partial q}{\partial x}(x_{n+1})\right]. \\
       \end{array}
\right.
\end{eqnarray}
In many applications, a special case is that $R(x)$ is a constant matrix, then the first equation is explicit since $\frac{\partial(p_{n+1/2}^TRp_{n+1/2})}{\partial x}(x_n)=0$. That is,
\begin{eqnarray}
p_{n+1/2}=\left[I_d+\frac{h}{2}\left[\frac{\partial f}{\partial x}(x_n)\right]^T\right]^{-1}\left[p_n-\frac{h}{2}\frac{\partial q}{\partial x}(x_n)\right].
\end{eqnarray}
Note that $I_d+\frac{h}{2}(\frac{\partial f}{\partial x}(x_n))^T$ is invertible since $h$ is small. Hence in \eqref{e:sv-s}, the second equation is the only implicit equation.

Symplectic algorithms have favourable long term behaviours such as energy conservation.
Assume that a Hamiltonian $H:D\to \mathbb R$ ($D\subset \mathbb R^{2d}$) is analytic. Suppose that $\Phi_h(y)$ is the St\"ormer-Verlet method with step size $h>0$. If the numerical solution stays in some compact set $K\subset D$, then there exists $h_0$ such that
\begin{eqnarray}\label{e:stable-sv}
H(x_n,p_n)=H(x_0,p_0)+O(h^2),
\end{eqnarray}
in exponential large interval $0<nh\le e^{h_0/(2h)}$.
For example, the Hamiltonian of the example in Section \ref{s:example} is analytic. Moreover, in concrete problem the constant $h_0$ can be computed explicitly \cite[Section 8.1]{hairer2006geometric}. For more general symplectic algorithms of various orders, similar energy estimates hold. See e.g. \cite[Section 8.1]{hairer2006geometric}.
As pointed out in \cite{sakamoto2008, sakamoto2013case}, the value of the Hamiltonian (i.e. the energy) is usually used as a measure for the accuracy of the approximate solutions. That is, if $H(x_n,p_n)$ is not close to $H(x_0,p_0)$, then the solution trajectory leaves the stable manifold. The estimates \eqref{e:stable-sv} indicate that the steps of symplectic algorithm can be controlled by the value of Hamiltonian well. In practice, we set $|H(x_n,p_n)-H(x_0,p_0)|$ along the numerical trajectories to satisfy certain accuracy $\delta>0$. Once $|H(x_n,p_n)-H(x_0,p_0)|>\delta$, the numerical computation will stop and record the time. Such technique is called \emph{Hamiltonian check} as in \cite{sakamoto2013case}.

\subsection{Computation procedure}\label{s:computation}
We are now in the position to summarize the computation procedure.

\emph{Step 1. Transform \eqref{e:ham-flow} into a system of form \eqref{e:Ham-transform}.} To apply the iterative method, we transform the Hamiltonian system into the separated form \eqref{e:Ham-transform} by a coordinates transformation as in Lemma \ref{l:diag}.

\emph{Step 2. Compute the local approximate stable manifold by iteration}. We give a precise estimate of the radius $\rho$ of $\xi$ which makes the sequences $\{x_k\}$ and $\{p_k\}$ convergent by Theorem \ref{t:iteration}.
Using numerical methods (e.g., Runge-Kutta method), \eqref{e:iterative} is solved for different points $\xi\in\mathbb S_{\rho}$. Here the number of points $\xi$ can be properly chosen in concrete problems. Then we get a local approximate stable manifold $\Lambda_k=\{(x_k(t,\xi),p_k(t,\xi))\,|\,t\ge 0,\,\xi\in \mathbb S_\rho\}$. Note that the error can be controlled by $k$ and $\rho$ by Theorem \ref{t:iteration}.

\emph{Step 3. Extend the local approximate stable manifold by symplectic algorithm.} Rewrite $\Lambda_k$ in the original coordinates by $\hat\Lambda_k=T\Lambda_k$ where $T$ is the coordinates transform given by Lemma \ref{l:diag}. Taking advantage of symplectic algorithm such as the St\"ormer-Verlet method, symplectic Runge-Kutta method of various orders, we solve the initial value problem \eqref{e:after-ct-extension-app-H}.
Then we find a larger approximate stable manifold. We shall use the Hamiltonian check to indicate that the trajectories stay close to the exact stable manifold.

\begin{remark}
In practice, we can find an applicable approximate stable manifold as follows: Firstly, a proper set of $\xi\in \mathbb S_{\rho}$ is chosen according to the concrete problem. From Theorem \ref{t:iteration}, a corresponding set of local curves in local approximate stable manifold is obtained. Then, extending the local curves by solving \eqref{e:after-ct-extension-app-H} numerically, we find a set of longer curves in the global approximate stable manifold.
Using $d$-dimensional interpolation method or $d$-variable polynomial fitting (see e.g. \cite{sakamoto2013case}), we obtain an approximate function $\tilde p(x)$ whose graph is an approximate stable manifold. For a detailed example, see Section \ref{s:control} below.
\end{remark}
\begin{remark}\label{r:sak-van}
In \cite{sakamoto2008}, \cite{sakamoto2013case}, the local approximate stable manifold is extended by using negative $t$ in \eqref{e:iterative} (or, equivalently, \eqref{e:iterative-int}) and taking more iterative times $k$. In our approach, we first generate a local approximate stable manifold by \eqref{e:iterative} with $t\ge 0$, then extend this local approximate stable manifold by solving the initial value problem \eqref{e:after-ct-extension-app-H} for negative $t$. We should emphasize that in our approach, we do not use negative $t$ in the iterative procedure \eqref{e:iterative}. This can avoid the divergent case of iterative sequence for negative $t$ as in \cite{sakamoto2008}, \cite{sakamoto2013case} and also reduce the computation cost.
\end{remark}

\section{Example}\label{s:example}
In this section, we apply the symplectic algorithm to a 2-dimensional optimal control problem with exponential nonlinearity. The existed methods may be difficult to applied to the systems with such kind of complicated nonlinearities as showed in \cite{sakamoto2013case,horibe2017optimal, horibe2018nonlinear}.

Throughout this section, we use the following notations.
$k$: the iterative times for local approximate stable manifold as in \eqref{e:iterative};
$\xi$: the initial condition given in Theorem \ref{t:iteration};
$h_-$: the step size of the numerical methods for $t<0$ by \eqref{e:after-ct-extension-app-H}.
RK45 method is from \emph{scipy.integrate.RK45} in python environment. This algorithm is based on error control as 4-order Runge-Kutta method, and its steps are taken using the 5-order accurate formula. For more details this algorithm, see \cite{dormand1980family, shampine1986some}. We should mention that other numerical environment such as Matlab also includes a similar algorithm named by ode45. In what follows, we implement RK45 in iterative procedure \eqref{e:iterative} for the local stable manifold.

Let us consider the following optimal control problem with exponential nonlinearity:
\begin{eqnarray}\label{e:2-dim}
\left\{\begin{array}{l}
         \dot x_1=e^{x_2}-1+u_1, \\
         \dot x_2=-(x_1+\frac{1}{3}x_1^3)+u_2,
       \end{array}
\right.
\end{eqnarray}
where $u=(u_1,u_2)^T$ is the feedback control function.
Let
$$
f(x)=\left[
         e^{x_2}-1,
         -(x_1+\frac{1}{3}x_1^3)
     \right]^T
.$$
Define the cost function by $$\int_0^{\infty}\frac{1}{2}[x^T(t)x(t)+u^T(t)u(t)]dt.$$ Then the corresponding Hamilton-Jacobi equation is
$$
H(x,p)=p^Tf(x)-\frac{1}{2}p^TRp+\frac{1}{2}x^TQx=0,
$$
where $p=\nabla V$ is the gradient of the value function in column form, $x=(x_1,x_2)^T$, $R=I_2$,
$Q=I_2$. Then by \eqref{e:u*12} the optimal controller is $u=-p$.
The associated Hamiltonian system is
\begin{eqnarray}\label{e:Ham-2-dim}
\left\{\begin{array}{l}
                   \dot x = f(x)-R p \\
                   \dot p = -\left[\frac{\partial f}{\partial x}\right]^Tp-Qx,
                 \end{array}\right.
\end{eqnarray}
where $\left[\frac{\partial f}{\partial x}(x)\right]^T=\left[
                                                           \begin{array}{cc}
                                                             0 & -1-x_1^2 \\
                                                             e^{x_2} & 0 \\
                                                           \end{array}
                                                         \right]
$.
Then $(0,0)$ is an equilibrium and the Hamiltonian matrix is given by
$\left[
                   \begin{array}{cc}
                     A & -R \\
                     -Q & -A^T \\
                   \end{array}
                 \right],
$
where $A= \left[
            \begin{array}{cc}
              0 & 1 \\
              -1 & 0 \\
            \end{array}
          \right].
$ The stabilizing solution $P$ of the Riccati equation is $I_2$. From Lemma \ref{l:diag}, we find a matrix
$
T=\left[
            \begin{array}{cc}
              I_2 & -0.5I_2 \\
              I_2 & 0.5I_2 \\
            \end{array}
          \right]\notag
$
such that
the Hamiltonian system \eqref{e:Ham-2-dim} becomes a separated form
\begin{eqnarray}\label{e:Ham-coor-trans}
\left[
              \begin{array}{c}
                \dot{\bar x} \\
                \dot{\bar p} \\
              \end{array}
            \right]=\left[\begin{array}{cc}
                                B & 0 \\
                                0 & -B^T
                              \end{array}\right]
                              \left[
              \begin{array}{c}
                \bar x \\
                \bar p \\
              \end{array}
            \right]
            +\left[
  \begin{array}{c}
    n_s(\bar x,\bar p) \\
    n_u(\bar x,\bar p) \\
  \end{array}
\right],
\end{eqnarray}
where $P$ is the stabilizing solution of the Riccati equation \eqref{e:Riccati}, $B=A-RP$, $x=x(\bar x,\bar p), p=p(\bar x,\bar p)$ are defined by
$
\left[
  \begin{array}{c}
    x \\
    p \\
  \end{array}
\right]=T\left[
           \begin{array}{c}
             \bar x \\
             \bar p \\
           \end{array}
         \right],
$ and
$
\left[
  \begin{array}{c}
    n_s(\bar x,\bar p) \\
    n_u(\bar x,\bar p) \\
  \end{array}
\right]:=T^{-1}\left[
\begin{array}{c}
f(x)-Ax \\
-\left[\frac{\partial f}{\partial x}\right]^Tp+A^Tp \\
\end{array}
\right].$
Here the eigenvalues of $B$ are $-1\pm i$.
Then we consider an iterative procedure as \eqref{e:iterative}
with $\bar x_0 = e^{B t}\xi, \,\bar p_0 = 0$, and
$
\bar x_k(0)=\xi,\,\bar  p_k(+\infty)=0
$ for $k=1,2,3,\cdots$.

For \eqref{e:Ham-coor-trans} we can choose
$M(L)=(9/4)L$ for $L>2/3$ and $M(L)=3/2$ for $L\le 2/3$. Moreover, $a=1, \,b=1$. Hence by Theorem \ref{t:iteration}, if $M(L)L\ge 3/8$, then $|\bar x_k|+|\bar p_k|<L$ for all $k=1,2,\cdots$. That is, we can choose $L\ge 1/4$. Hence for $|\xi|\le \frac{3}{16M(L)}\approx 0.125$, $\{\bar x_k(t, \xi)\}$ and $\{\bar p_k(t, \xi)\}$ converge to exact solution $\bar x(t, \xi)$ and $\bar p(t, \xi)$ respectively. In what following, we shall take $\xi$ in the sphere $\mathbb S_{0.12}$.  That yields a local approximate stable manifold $\Lambda_k=\{(\bar x_k(t,\xi), \bar p_k(t,\xi))\,|\,t\ge 0,\,\xi\in \mathbb S_{0.12}\}$ and $\partial \Lambda_k=\{(\xi,\bar p_k(0,\xi))\,|\,\xi\in \mathbb S_{0.12}\}$.

Next, as in Section \ref{s:computation}, the local stable manifold $\Lambda_k$ will be enlarged to a global one $\Lambda_{k,g}$ from \eqref{e:after-ct-extension-app-H}. We solve the initial value problem \eqref{e:after-ct-extension-app-H} by the St\"ormer-Verlet scheme \eqref{e:sv-s}. In our problem \eqref{e:Ham-2-dim}, algorithm \eqref{e:sv-s} becomes
\begin{eqnarray}
\left\{\begin{array}{l}
         p_{n+1/2}=\left[I_2+\frac{h}{2}\left[\frac{\partial f}{\partial x}(x_n)\right]^T\right]^{-1}\left[p_n-hx_n\right]\\
         x_{n+1} = x_n+\frac{h}{2}\left[f(x_n)+f(x_{n+1})-2p_{n+1/2}\right] \\
         p_{n+1}=p_{n+1/2}-\frac{h}{2}\left[\left[\frac{\partial f}{\partial x}(x_{n+1})\right]^Tp_{n+1/2}+x_{n+1})\right]. \\
       \end{array}
\right.\notag
\end{eqnarray}
Here the second equation is the only implicit equation which can be solved by Newton's iteration method.

\subsection{Comparison with other methods}
\begin{figure}[htbp]
\begin{center}
\subfigure{
\includegraphics[width=0.48\textwidth]{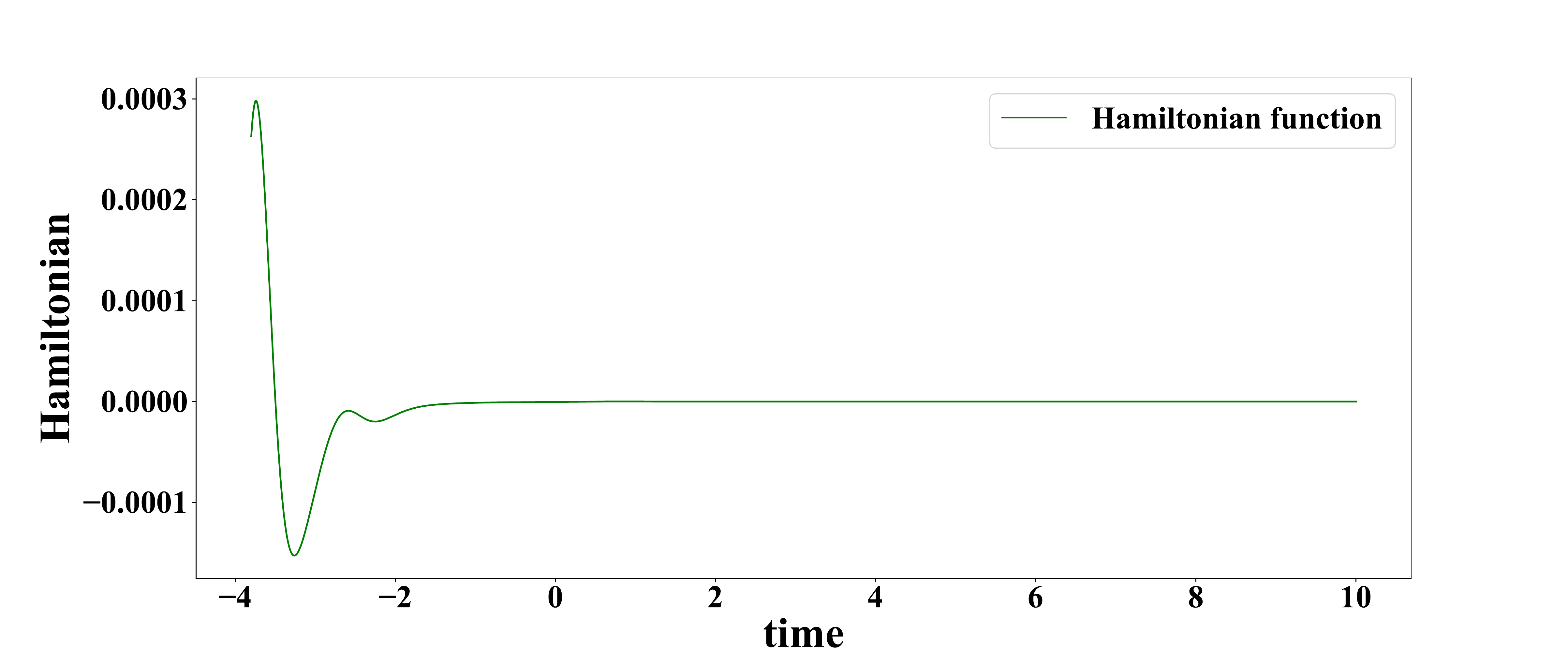}}
\subfigure{
\includegraphics[width=0.48\textwidth]{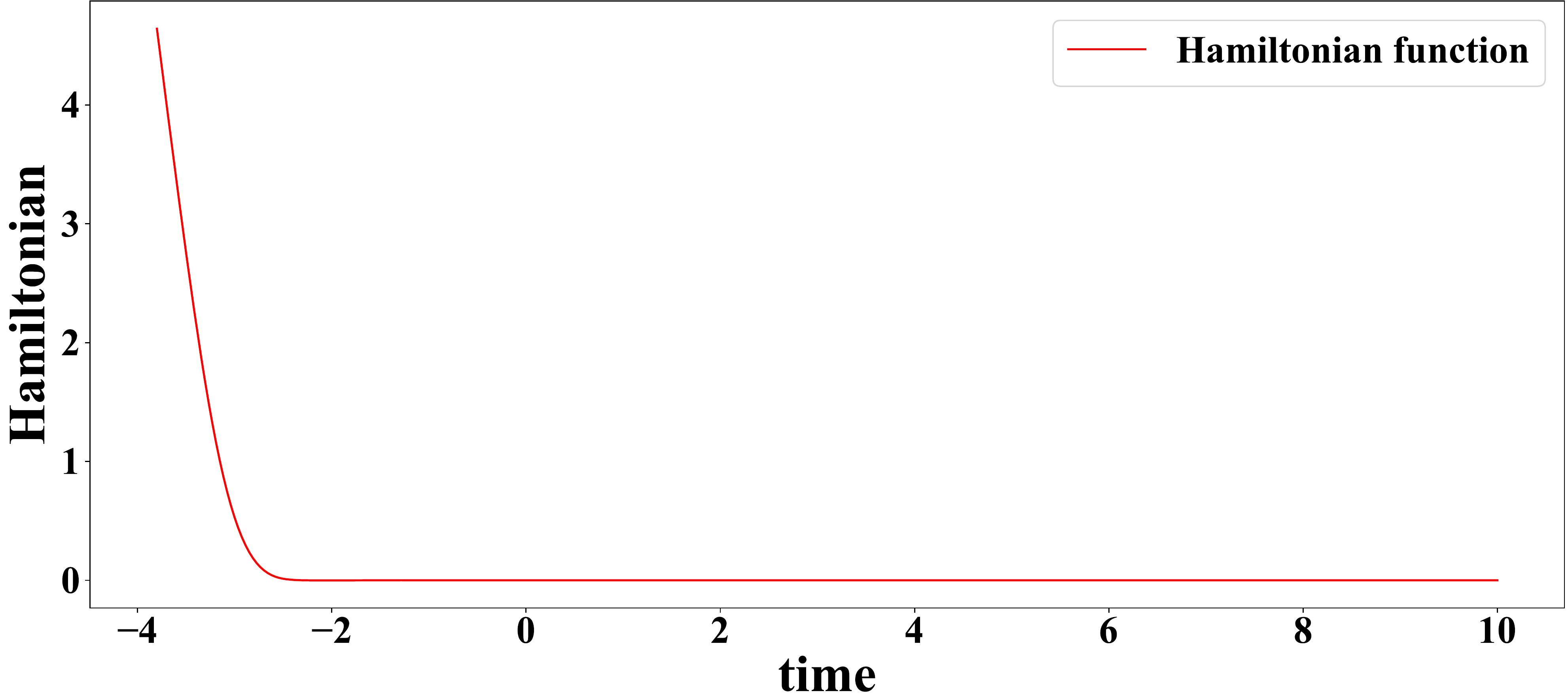}}
\subfigure{
\includegraphics[width=0.48\textwidth]{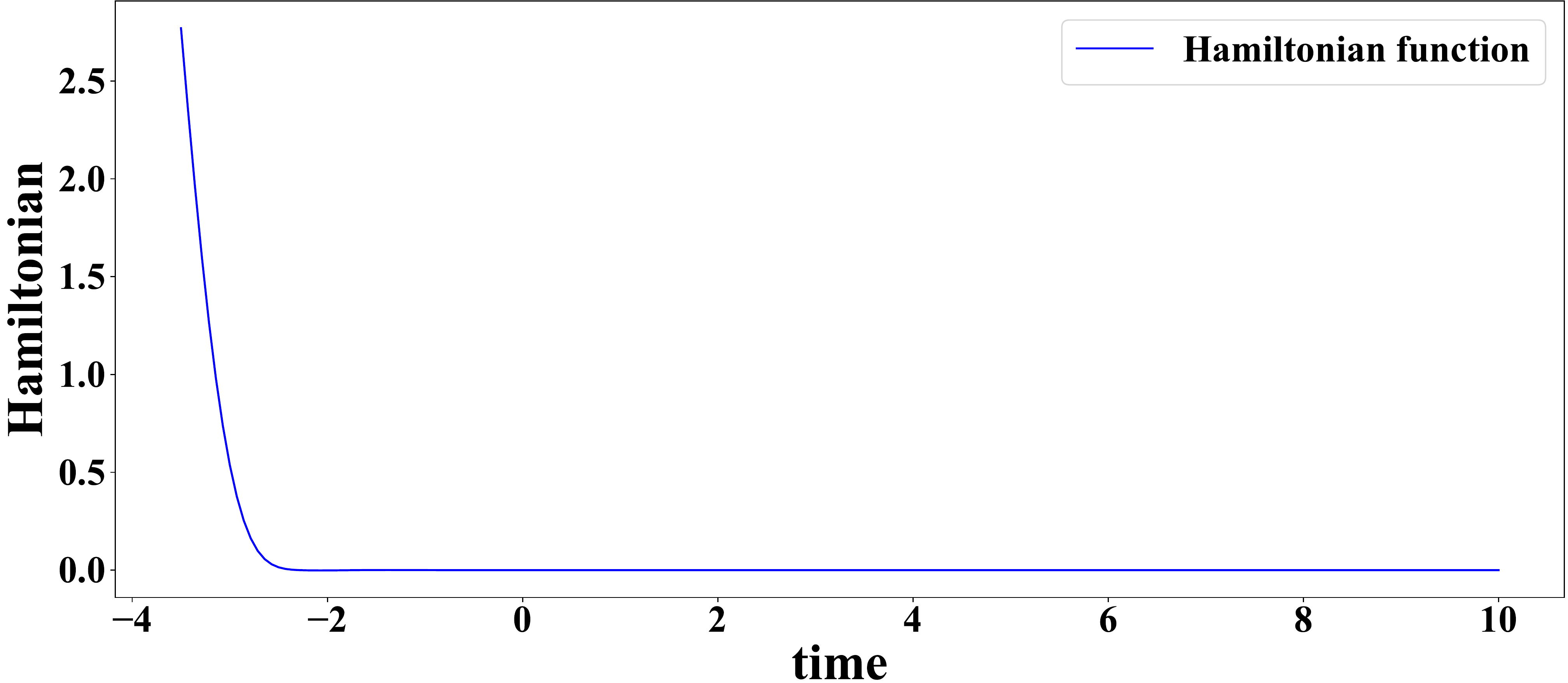}}
\end{center}
\hfill
\vspace{-0.6cm}
\caption{ \footnotesize Hamiltonian values along extended curves with $\xi=(0.12\times\sqrt{2}/2,0.12\times\sqrt{2}/2)$. The first subfigure is the St\"ormer-Verlet scheme with $k=3$, $h_-=-0.005$, $t\in [-3.8,10]$; the second one corresponds to the RK45 extension method with $k=3$, $t\in [-3.8,10]$; the third one is from the extension method of \cite{sakamoto2008, sakamoto2013case} with $k=50$, $t\in [-3.5,10]$. We should point out that in the third subfigure, the time interval is smaller since the Hamiltonian value from the extension method of \cite{sakamoto2008, sakamoto2013case} blows up to infinite on $[-3.8,10]$.}
\label{f:exp-1}
\vspace{-0.8cm}
\end{figure}

To show the effectiveness of our approach, the extension result of the St\"ormer-Verlet scheme is compared with other two methods, i.e., the first one is solving extension initial problem \eqref{e:after-ct-extension-app-H} by the RK45 algorithm, and the second one is the method of iteration for negative time in problem \eqref{e:iterative-int} as in Sakamoto-van der Schaft \cite{sakamoto2008, sakamoto2013case}. For the differences of these methods, see Remark \ref{r:sak-van} above. Figure \ref{f:exp-1} illustrates the values of the Hamiltonian function along the approximate curves from the three methods with $\xi=(0.12\times\sqrt{2}/2,0.12\times\sqrt{2}/2)\in\mathbb S_{0.12}$. It shows that the St\"ormer-Verlet method is much better than the RK45 extension method and the extension approach of Sakamoto-van der Schaft \cite{sakamoto2008, sakamoto2013case}.

\begin{figure}[htbp]
\begin{center}
\subfigure{
\includegraphics[width=0.35\textwidth]{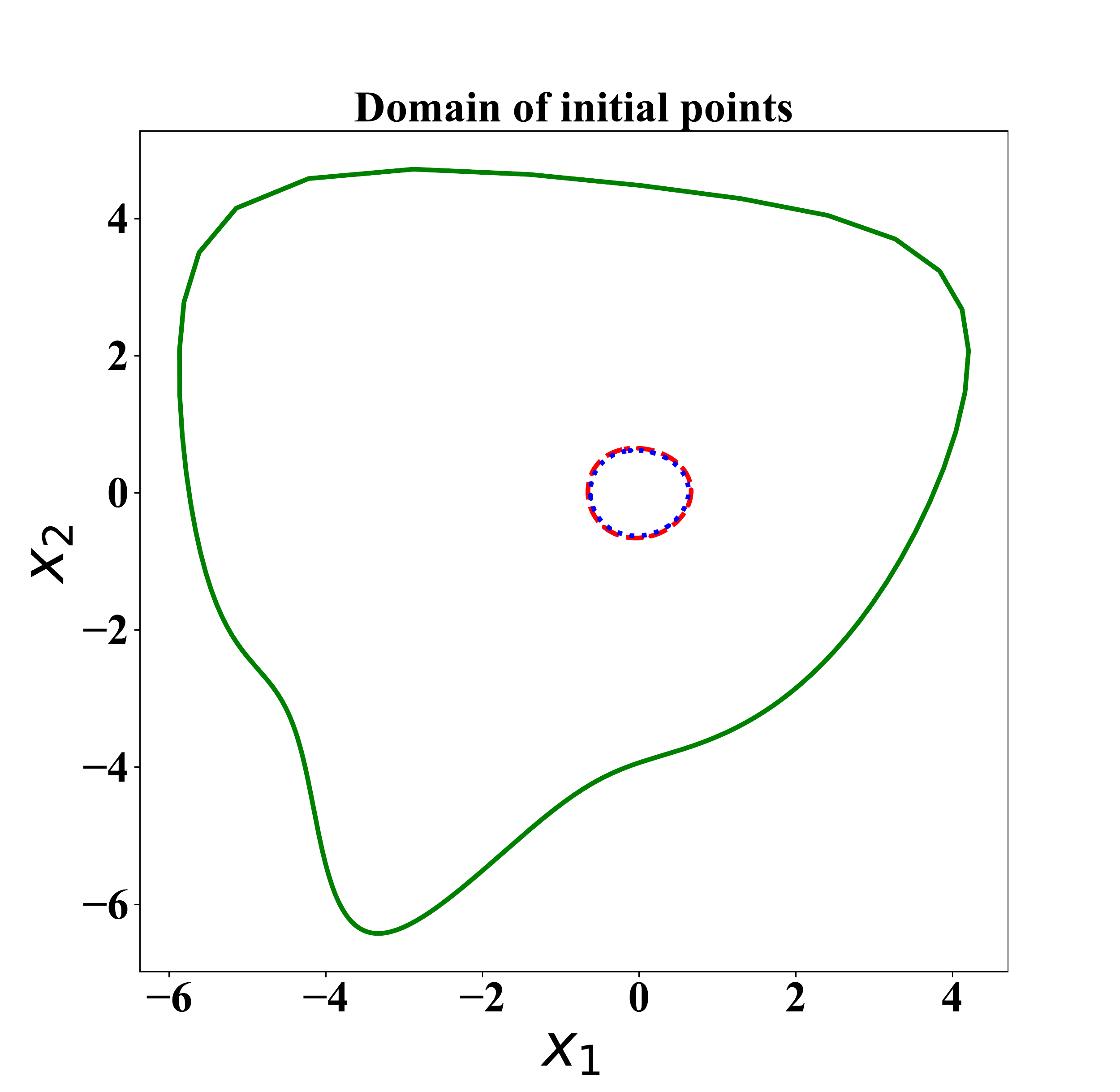}}
\end{center}
\hfill
\vspace{-0.8cm}
\caption{ \footnotesize The domains generated by projections of extended approximate local stable manifold to $x$-plane. The largest domain $\Omega$ (with green line boundary) is from the extension by the St\"ormer-Verlet scheme for $k=3$, $h_-=-10^{-3}$ and the absolute value of Hamiltonian less than $10^{-4}$. The domain with red dashed line boundary (resp., blue dotted line boundary) corresponds to the RK45 extension method for $k=3$ (resp., the extension method of \cite{sakamoto2008, sakamoto2013case} for $k=50$) with the absolute value of Hamiltonian less than $10^{-3}$.}
 \label{f:exp-domain}
 \vspace{-0.5cm}
\end{figure}

From \eqref{e:app-stable}, the approximate stable manifold $\Lambda_{k,g}$ is made up by all curves of solutions $(x_k(t,\xi),p_k(t,\xi))$, $\xi\in \mathbb S_{0.12}$. To construct the optimal controller, we need to find numerical solutions $(x_k(t,\xi),p_k(t,\xi))$ with Hamiltonian value under certain error tolerance. Projection of the approximate stable manifold (i.e. all extended curves) under certain error tolerance to $x$-plane is a domain. In Figure \ref{f:exp-domain}, we illustrate the domains of projection of all extended curves with $\xi\in \mathbb S_{0.12}$ to $x$-plane by the St\"ormer-Verlet method, RK45 and the method of Sakamoto-van der Schaft (\cite{sakamoto2008, sakamoto2013case}) under certain error tolerances. It is clear that the St\"ormer-Verlet method generates a much larger domain comparing with the other two extension methods.

We should point out that, as shown above, the method of Sakamoto-van der Schaft (\cite{sakamoto2008, sakamoto2013case}) requires much more iterative times, hence the computation cost is more expensive.
\begin{remark}
We should mention that the St\"ormer-Verlet scheme is just a two order symplectic algorithm. If we need a larger domain, or a smaller Hamiltonian value tolerance is required, then other symplectic algorithms such as symplectic Runge-Kutta schemes of higher orders (see e.g. \cite{hairer2006geometric}) can be implemented.
\end{remark}

\subsection{Construct optimal controller by the approximate stable manifold}\label{s:control}
\begin{figure}[htbp]
\begin{center}
\subfigure{
\includegraphics[width=0.48\textwidth]{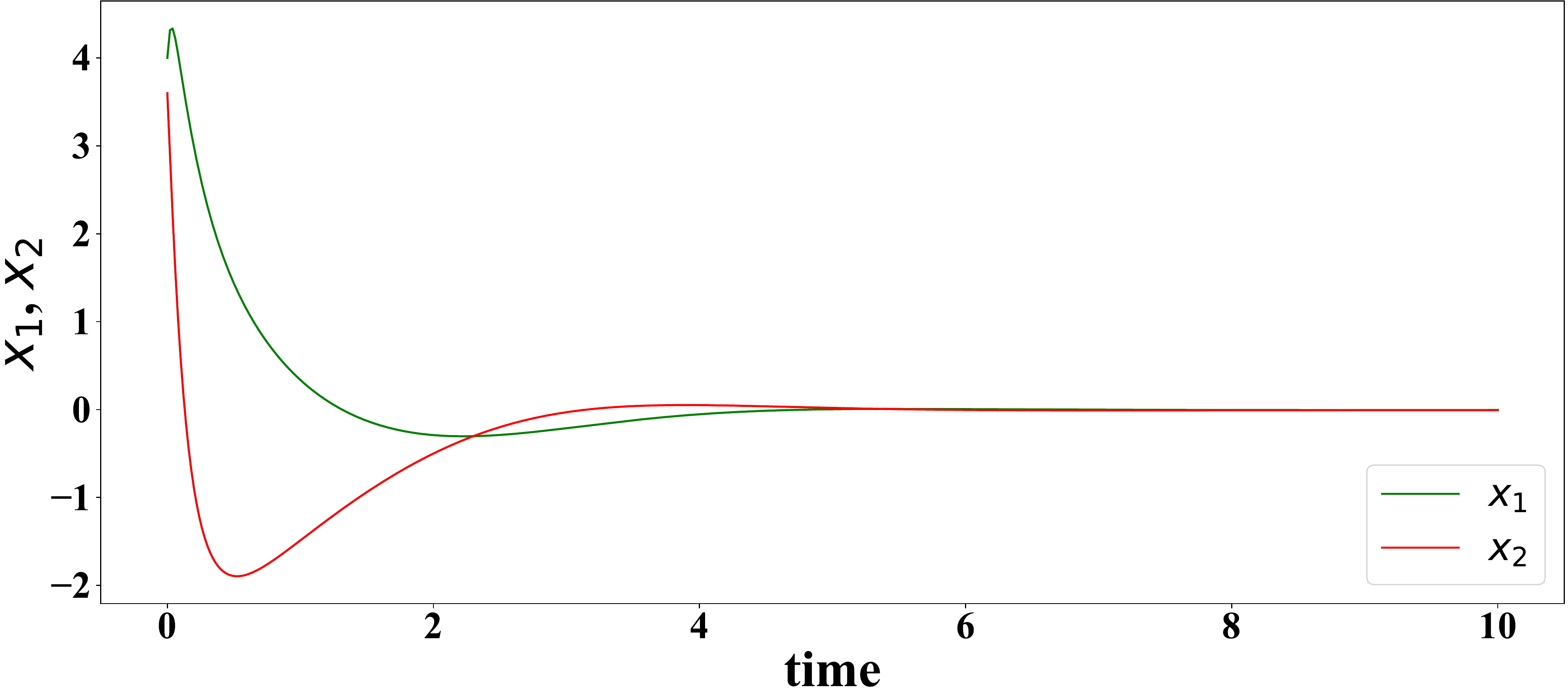}}
\subfigure{
\includegraphics[width=0.48\textwidth]{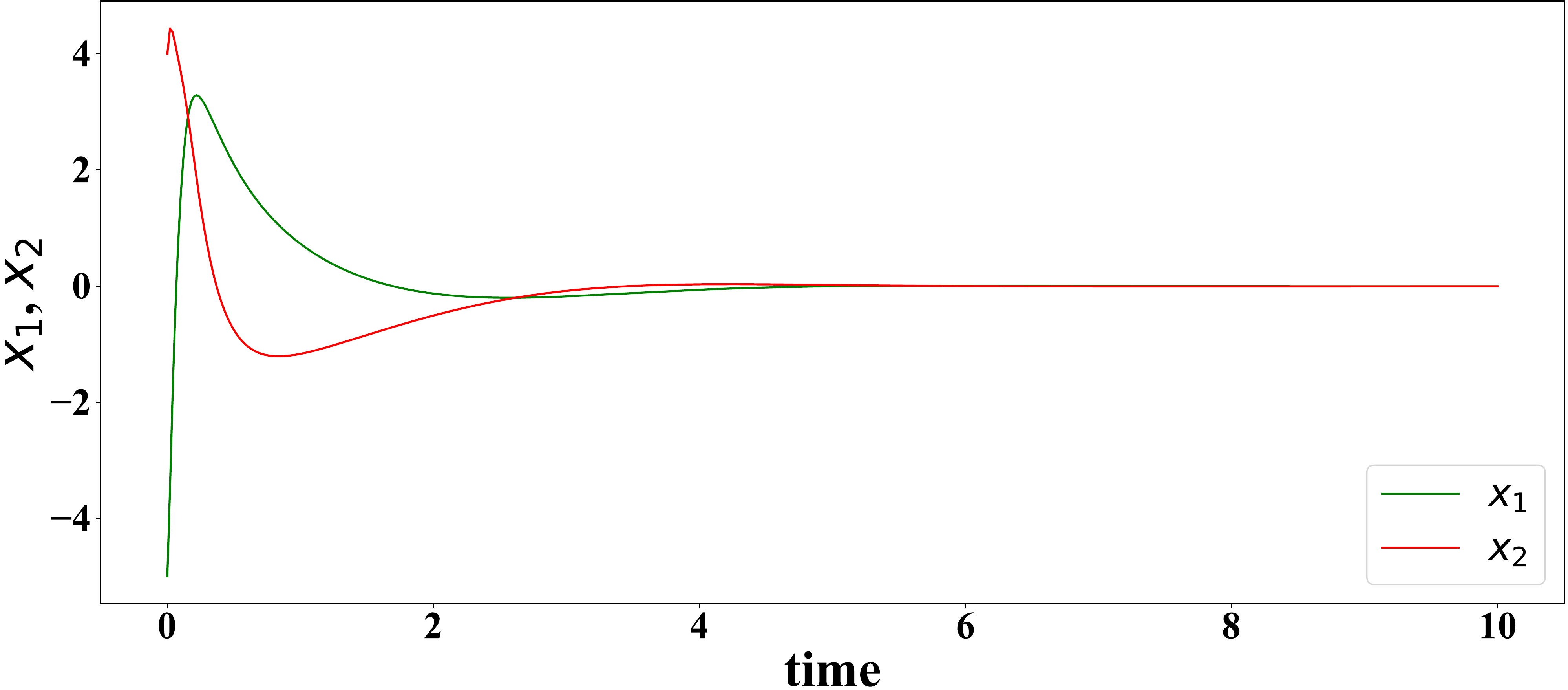}}
\end{center}
\hfill
\vspace{-0.8cm}
\caption{\footnotesize  The controlled trajectories for the initial point $(4,3.6)$ and $(-5,4)$.}
 \label{f:exp-controller}
\end{figure}

To find the approximate optimal controller $u$ by the stable manifold, we use the polynomial fitting method as in \cite{sakamoto2013case}. To be more precise, we first choose $200$ $\xi$s according to the uniform distribution in $\mathbb S_{0.12}$. Using the computation procedure in Section \ref{s:computation} based on St\"ormer-Verlet scheme for time interval $[-3.5, 10]$, $h_-=-10^{-3}$ and $k=3$, we find $200$ extended curves in the stable manifold with the Hamiltonian less than $10^{-4}$. The projections of these curves to $x$-plane belong to the domain $\Omega\subset \mathbb R^2$ (the largest domain in Figure \ref{f:exp-domain}).
Then we select 10 points (including the point at time $-3.5$ and the other 9 points chosen randomly in $(-3.5,0]$ according to the uniform distribution) on each curve and the point $(x,p)=(0,0)$. Hence we obtained a set $\mathcal W$ of sample points on the stable manifold. To approximate stable manifold, we assume a polynomial fitting of form
$$p_{\rm pol}(x)=\left[\sum_{0\le i,j\le 5}C^1_{ij}x_1^ix_2^j, \sum_{0\le i,j\le 5}C^2_{ij}x_1^ix_2^j\right]^T,$$
where the coefficients $C_{ij}^k$ are chosen by the least square technique based on the sample set $\mathcal W$. That is, $\{C_{ij}^k\,|\,i,j=1,\cdots,5; k=1,2\}$ are the coefficients to minimize
$$\sum_{(x_l,p_l)\in \mathcal W}(p_{\rm pol}(x_l)-p_l)^2.$$
Then, using \eqref{e:u*12}, the approximate optimal controller is given by $\tilde u(x)=-p_{\rm pol}(x)$.
Hence by system \eqref{e:2-dim}, the closed loop can be computed for any $x\in \Omega$. Here we demonstrate the controlled trajectories from the controller $\tilde u$ at some points in Figure \ref{f:exp-controller}. We can clearly see that the controller asymptotically stabilizes these trajectories.

\section{Conclusion}
In this note, using a symplectic algorithm for the stable manifolds of the Hamilton-Jacobi equations from infinite horizon optimal control, we construct a sequence of local approximate stable manifolds for Hamiltonian system at some hyperbolic equilibrium.
In our approach, we first construct an iterative sequence of local approximate stable manifolds at equilibrium based on a precise estimates for the radius of convergence.
Then, using symplectic algorithms, we enlarge the local approximate stable manifolds by solving initial value problem of the associated Hamiltonian system for $t<0$. This avoids the divergent case of the iterative sequence as in \cite{sakamoto2008} for unstable direction $t<0$, and the computation cost can be reduced.
Moreover, the symplectic algorithms have better long-time behaviours than general purpose numerical methods such as Runge-Kutta. The effectiveness of our algorithm is demostrated by an optimal control problem with exponential nonlinearity. Lastly, we should mention that the symplectic algorithm has potentials to be applied to more general Hamilton-Jacobi equations in control theory since such kind of equations have associated Hamiltonian systems as the characteristic systems.

\appendix

\subsection{Proof of Theorem \ref{t:iteration}}
Firstly, we prove that for $|\xi|\le \frac{3b}{16a^2M}$,
\begin{eqnarray}\label{e:alpha-beta-bounds}
|x_k(t)|\le \underline{\alpha} e^{-bt},\quad |y_k(t)|\le \underline{\beta} e^{-2bt}, \quad \forall t\ge 0,
\end{eqnarray}
where
\begin{eqnarray}\label{e:al-beta}
\underline{\alpha}=\frac{\frac{3a^2}{16}|\xi|^2}{g+\sqrt{g^2-\frac{a^2}{16}|\xi|^2}}+a|\xi|,\,
\underline{\beta}=\frac{\frac{a^2}{16}|\xi|^2}{g+\sqrt{g^2-\frac{a^2}{16}|\xi|^2}}.
\end{eqnarray}
Here $g=\frac{3b}{32aM}-\frac{a}{4}|\xi|\ge \frac{3b}{64aM}$.

In fact, By Assumption 1,
$|x_0(t)|\le a|\xi|e^{-bt}, |y_0(t)|=0$, $t\ge 0.$
Hence $\alpha_0 = a|\xi|$, $\beta_0=0$. For $k=1,2,\cdots$, we prove by inductive method. Assume that the claim holds for $k$. Then from \eqref{e:n-assumption} and \eqref{e:alpha-beta-bounds},
\begin{eqnarray}
|x_{k+1}(t)|
\le a|\xi|e^{-bt}+aMe^{-bt}\int_0^t e^{bs}(|x_k(s)+y_k(s)|^2)ds\notag\\
\le \left[\frac{aM}{b}(\alpha_k+\beta_k)^2+a|\xi|\right]e^{-bt},\quad \mbox{ and }\notag
\end{eqnarray}
\begin{eqnarray}
|y_{k+1}(t)|
&\le&aMe^{bt}\int_t^\infty e^{-bs}(|x_k(s)+y_k(s)|^2)ds\notag\\
&\le& \frac{aM}{3b}(\alpha_k+\beta_k)^2e^{-2bt}.\notag
\end{eqnarray}
Let $c_0=\frac{aM}{3b}$. Define $ \alpha_0=a|\xi|,\,\beta_0 = 0$, $\alpha_{k+1} = 3c_0(\alpha_k+\beta_k)^2+a|\xi|, \,\beta_{k+1} =c_0(\alpha_k+\beta_k)^2,$ for $k=1,2,3,\cdots$.
We should point out that the definition of $\alpha_{k}, \beta_{k}$ is different from that in \cite{sakamoto2008}.
Note first that $\alpha_1>\alpha_0$ and $\beta_1>\beta_0$. By mathematical induction, we can prove that $\alpha_{k+1}>\alpha_k$ and $\beta_{k+1}>\beta_k$. Next solving
\begin{eqnarray}
\left\{\begin{array}{l}
         \underline \alpha=3c_0(\underline \alpha+\underline \beta)^2+a|\xi|, \\
         \underline \beta=c_0(\underline \alpha+\underline \beta)^2,
       \end{array}
\right.
\end{eqnarray}
we have solution \eqref{e:al-beta}. Then \eqref{e:alpha-beta-bounds} follows.

Secondly, by a similar argument as in \cite{sakamoto2008}, we find that
\begin{eqnarray}\label{e:xk-yk}
|x_{k+1}(t)-x_k(t)|\le \gamma_k e^{-bt},\, |y_{k+1}(t)-y_k(t)|\le \varepsilon_k e^{-2bt}
\end{eqnarray}
where $\{\gamma_k\}$ and $\{\varepsilon_k\}$ satisfy $\gamma_1=\frac{a^3M|\xi|^2}{b} ,\,\varepsilon_1 = \frac{a^3M|\xi|^2}{3b}$ and $\gamma_{k+1}=\frac{a(\underline \alpha+\underline \beta)M}{b}(\gamma_k+\varepsilon_k)$, $\varepsilon_{k+1}=\frac{a(\underline \alpha+\underline \beta)M}{3b}(\gamma_k+\varepsilon_k)$.
Moreover, $\{\gamma_k\}$ and $\{\varepsilon_k\}$ are decreasing and $\lim_{k\to\infty}\gamma_k = 0$, $\lim_{k\to\infty}\varepsilon_k = 0$. Consequently, it holds that
\begin{eqnarray}\label{e:gamma+epsilon}
\gamma_k + \varepsilon_k &\le& \left[\frac{4}{3}\frac{a(\underline \alpha+\underline \beta)M}{b}\right]^{k-1} \frac{4a^3M|\xi|^2}{3b}.
\end{eqnarray}

Thirdly, we prove that if $M(L)L>\frac{3b}{8a}$, it holds that for all $k\in \mathbb N$,
\begin{eqnarray}\label{e:x_k-y_k-bounds}
|x_k(t)|+|y_k(t)|\le L,\quad \forall t\in [0,\infty).
\end{eqnarray}

In fact,
from \eqref{e:al-beta} and $|\xi|\le\frac{3b}{16a^2M}$, it holds that
$
\underline\alpha\le \frac{21}{64}\frac{b}{aM},\quad\underline\beta\le \frac{3}{64}\frac{b}{ aM}.
$
Hence
$
|x_k(t)|\le \frac{21}{64}\frac{b}{aM} e^{-bt}.
$
Similarly,
$
|y_k(t)|\le \frac{3}{64}\frac{b}{ aM}e^{-2bt}.
$
Therefore, we obtain that
$
|x_k(t)|+|y_k(t)|\le \frac{3}{8}\frac{b}{aM} e^{-bt}.
$
Then since $M(L)L>\frac{3b}{8a}$, \eqref{e:x_k-y_k-bounds} holds. Since $M=M(L)$ is increasing with respect to $L$, the conclusion holds.

Finally, let $k\in \mathbb N$ be any fixed number. From \eqref{e:xk-yk} and \eqref{e:gamma+epsilon},
for all $j\in \mathbb N$, it holds that
\begin{eqnarray}\label{e:x-y-j}
&&|x_{k+j}(t)-x_k(t)|\\
&\le& \left[\frac{4}{3}\frac{a(\underline \alpha+\underline \beta)M}{b}\right]^{k-1}\frac{4a^3M|\xi|^2}{3(b-a(\underline \alpha+\underline \beta)M)} e^{-bt},\notag\\
&&|y_{k+j}(t)-y_k(t)|\notag\\
&\le& \left[\frac{4}{3}\frac{a(\underline \alpha+\underline \beta)M}{b}\right]^{k-1}\frac{4a^3M|\xi|^2}{3b-a(\underline \alpha+\underline  \beta)M}e^{-2bt}.\notag
\end{eqnarray}
Here we used the fact that $\frac{4}{3}\frac{a(\underline \alpha+\underline \beta)M}{b}\le 1/2$. Rewriting \eqref{e:x-y-j} in form \eqref{e:x-y}, the conclusions of this theorem hold since $j$ is arbitrary. This completes the proof. \hfill $\Box$

\end{document}